\documentclass{article}
\usepackage{graphicx} % Required for inserting images
\usepackage[utf8]{inputenc}
\usepackage{authblk}
\usepackage{blindtext}
\usepackage{geometry}
\usepackage{amsmath,amssymb,amsthm}
\usepackage{graphicx,url,comment}
\usepackage{mathrsfs}
\usepackage{mathtools}
\usepackage{xcolor}
\usepackage[toc,page]{appendix}
\usepackage{tikz}
\usepackage{authblk}
\usepackage[colorlinks=true, urlcolor=blue, linkcolor=blue, citecolor=blue]{hyperref}

\theoremstyle{definition}
\newtheorem{theorem}{Theorem}[section]
\newtheorem{proposition}[theorem]{Proposition}
\newtheorem{definition}[theorem]{Definition}
\newtheorem{lemma}[theorem]{Lemma}
\newtheorem{question}[theorem]{Question}
\newtheorem{example}[theorem]{Example}
\newtheorem{remark}[theorem]{Remark}
\newtheorem{corollary}[theorem]{Corollary}

\newtheorem{exercise*}[exercise]{Exercise*}

\newcommand{\K}{\mathbb{K}}
\newcommand{\PP}{\mathbb{P}}
\newcommand{\F}{\mathbb{F}}
\newcommand{\Z}{\mathbb{Z}}
\newcommand{\Cc}{\mathcal{C}}
\newcommand{\Gg}{\mathcal{G}}
\newcommand{\Pp}{\mathcal{P}}
\newcommand{\Dd}{\mathcal{D}}

\newcommand{\reg}{\mathrm{reg}}
\newcommand{\nbl}{\mathrm{nb}}

\title{A combinatorial description of when a self-associated set of points fails to be arithmetically Gorenstein\thanks{The first and the second authors were partially supported by Grant PID2022-138906NB-C21 funded by MCIN/AEI/ 10.13039/501100011033 and  by ERDF/UE, and by Grant CLU-2025-1-02-IMUVA funded by Department of Education of the Junta de Castilla y Le\'on and FEDER Funds. The third author was supported by the Swiss National Science Foundation under the grants P500PT-222344 and P5R5-2\_239500}}
\author[1]{Gonzalo Rodríguez-Pajares \thanks{Email: gonzalo.rodriguez.pajares@uva.es}}
\author[1]{Diego Ruano \thanks{Email: diego.ruano@uva.es}}
\author[2]{Flavio Salizzoni \thanks{Email: flavio.salizzoni@math.uzh.ch}}
\affil[1]{IMUVa-Mathematics Research Institute, Universidad de Valladolid, Spain}
\affil[2]{Max Planck Institute for Mathematics in the Sciences, Leipzig, Germany}
\affil[2]{University of Zurich, Switzerland}
 
\date{}

\providecommand{\keywords}[1]
{
  \textbf{\textit{Key words: }} #1
}
\providecommand{\subclas}[1]
{
  \textbf{\textit{Mathematics Subject Classification: }} #1
}

\begin{document}
\maketitle
\begin{abstract}
We give a combinatorial proof that the set of points associated to a self-dual code with no proportional columns is arithmetically Gorenstein if and only if the code is indecomposable. This answers a question asked by Toh{\u{a}}neanu. We do so by providing a combinatorial way to compute the dimension of the Schur square of a self-dual code through a zero-one symmetrization of a systematic generator matrix. Our approach also allows us to compute the Gorenstein defect. As a consequence, we obtain a combinatorial characterization of arithmetically Gorenstein self-associated sets of points over an algebraically closed field. We also show that, for any finite field, almost all self-dual codes of large length are indecomposable, and hence their associated sets of points are arithmetically Gorenstein.

\,

\noindent\keywords{Self-associated set, Arithmetically Gorenstein, Self-dual codes, Schur product, Gale transform.}
%11T71

\noindent\subclas{94B05, 13H10, 14N10.}
\end{abstract}
 
\begingroup
\renewcommand\thefootnote{}
%\footnotetext{}
\endgroup

\section{Introduction}
The Gale transform is an involution that maps an ordered set $\Pi$ of $n$ points in $\mathbb{P}^{k-1}$ to an ordered set $\Pi'$ of $n$ points in $\mathbb{P}^{n-k-1}$, defined up to projective linear transformations. An ordered set of points is called self-associated or self-dual if it is invariant under the Gale transform; such a set necessarily consists of $2k$ points in $\mathbb{P}^{k-1}$. Although the notion of Gale transform as we know it today was given by Coble in~\cite{coble1922associated}, the history of self-associated points is longer and can be traced back to Pascal's early studies on points on conics~\cite{taton1955essay}. This concept was then formalized by Castelnuovo in~\cite{castelnuovo1889certi} a couple of centuries later. It was characterized in~\cite[Theorem 7.3]{eisenbud2000projective} that, over an algebraically closed field, a nondegenerate self-associated set of points is arithmetically Gorenstein if and only if it fails exactly by one to impose independent conditions on quadrics. The importance of studying arithmetically Gorenstein sets of points stems from the connection with free resolutions and Serre duality~\cite{isenbud1999resolutions}.

This theory can be reformulated in terms of linear codes. The Gale transform corresponds to taking the dual code, and the columns of a generator matrix of a self-dual code form a self-associated set of points. Conversely, a self-associated set of points gives a code which is equivalent to a self-dual code when the field is algebraically closed, but not in general, see Remark~\ref{remark:algcl}.

In~\cite{tohuaneanu2024commutative} S.O.I. Toh{\u{a}}neanu asked whether the points given by the columns of a generator matrix of a self-dual code, with no proportional columns, always form an arithmetically Gorenstein set. This is particularly interesting in coding theory, as in this case one has a lower bound for the minimum distance in terms of the Castelnuovo–Mumford regularity of the coordinate ring of the set of points, see \cite{tohuaneanu2024commutative} for further details. In this work, we answer this question by providing a counterexample and by giving a combinatorial characterization of when a self-associated set of points is arithmetically Gorenstein, along with a description of the Gorenstein defect: the set of points is arithmetically Gorenstein if and only if the code is indecomposable, and in general the Gorenstein defect is the number of indecomposable components of the code minus one.

We learned after a first version of this work was written that this characterization already appeared in a scheme-theoretic framework in~\cite[Corollary 7.7]{eisenbud2000projective}. On the coding-theoretic side, part of the argument used to prove \cite[Corollary 7.7]{eisenbud2000projective} can also be obtained from a result that Randriambololona attributes to Couvreur and Tillich~\cite[Corollary 2.42]{randriam15course}, see Remark~\ref{remark:randriam1}. Moreover, shortly before our first version was posted, Kreuzer~\cite[Theorem 4.4]{kreuzer2026codeequivalencepointset} had already proved, by algebraic methods based on the canonical module, that an indecomposable iso-dual code gives rise to an arithmetically Gorenstein set of points; we were not aware of his work at that time. Our contribution is a different, elementary and self-contained proof, combinatorial in nature: we show that the dimension of the Schur square of a self-dual code is $2k$ minus the number of connected components of a bipartite graph attached to a systematic generator matrix of the code. This provides a clean connection between~\cite[Theorem 7.3]{eisenbud2000projective} and~\cite[Corollary 2.42]{randriam15course}, and it gives at the same time the Gorenstein defect and a simple procedure to compute it. Since the proof only uses linear algebra and  graph theory, we hope that it makes this connection accessible to readers from both coding theory and combinatorial commutative algebra.

The key ingredients in obtaining these results involve the Schur square of the linear code, that is, the componentwise product, and its relation with the Veronese map. Namely, for a self-dual code of dimension $k$ with no proportional columns, the columns of a generator matrix form an arithmetically Gorenstein set if and only if the dimension of the Schur square is $2k-1$. Thus, we will essentially compute the dimension of the Schur square of the linear code considered. Note that the computation of the Schur square has important applications in cryptography, in particular in code-based cryptography, and its relation with the Gale transform has already been examined, see \cite{randriam25syzygy} for further details. The combinatorial description relies on the number of connected components of a graph derived from a zero–one symmetrization of a systematic generator matrix. This formulation leads to a computation that is not computationally intensive: once the systematic generator matrix is known, the number of components can be found in $O(n^2)$ operations.

We remark that our results extend beyond the coding theory framework: although motivated by coding theory, they apply to arbitrary algebraically closed fields since the constructed graph and techniques are field independent. In particular, over an algebraically closed field we obtain a complete characterization of arithmetically Gorenstein self-associated sets of points: a self-associated set of $2k$ distinct points in $\mathbb{P}^{k-1}$ is arithmetically Gorenstein if and only if the linear code obtained by evaluating linear forms at the points is indecomposable (Corollary~\ref{corollary:alg closed}), i.e., it cannot be written as the direct sum of two nontrivial codes.

Finally, we connect our results with the theory of indecomposable linear codes \cite{BOGART1977333, Haluk1989phd}. We show that Gorenstein sets of points arising from self-dual codes over a finite field correspond to indecomposable self-dual linear codes. Moreover, an indecomposable self-dual code of length at least $4$ has no proportional columns. It is known that for sufficiently large lengths, almost all binary self-dual codes are indecomposable \cite{Haluk1989phd}. We extend this result to any finite field, and thus almost all self-dual codes give rise to arithmetically Gorenstein sets of points. In the case of self-dual codes with proportional columns, our computation can still be used to determine the dimension of the squared code. 

In Section \ref{Section 2}, we provide a self-contained introduction to all the concepts and results mentioned above, and we give a counterexample to the question posed in \cite{tohuaneanu2024commutative}. Section \ref{Section 3} contains the main results, where we characterize when the points of a self-dual code are arithmetically Gorenstein and give a combinatorial description of the Gorenstein defect. Finally, in Section \ref{Section 4}, we study the number of indecomposable self-dual codes and show that asymptotically almost all self-dual codes over a finite field are indecomposable. We illustrate all results with detailed examples.

\section{Preliminaries}\label{Section 2}

Let $\K$ be a field, $n$ and $k$ natural numbers, and $\PP^{k-1}(\K)$ the projective space of dimension $k-1$ over $\K$.  Given an ordered set of $n$ points $\Pi\subset\PP^{k-1}(\K)$ we say that a matrix $G\in\K^{k\times n}$ represents $\Pi$ if its columns are representatives in $\K^k$ of the points of $\Pi$. A set of points is said to be nondegenerate if it is not contained in a hyperplane. The Gale transform of $\Pi$ is defined as follows. 
\begin{definition}\label{def:Gale}
Let $\Pi\subset\PP^{k-1}(\K)$ and $\Pi'\subset\PP^{n-k-1}(\K)$ be ordered nondegenerate sets of $n$ points represented by the matrices $G\in\K^{k\times n}$ and $G'\in\K^{(n-k)\times n}$, respectively. We say that $\Pi'$ is the Gale transform of $\Pi$ if there exists a nonsingular diagonal matrix $D$ such that $GDG'^T=0$. An ordered set of points is called self-associated or self-dual if it is invariant under the Gale transform.
\end{definition}

Let $\Cc$ be an $[n,k]$ linear code over $\K$. Note that a code is said to be projective if its dual has minimum distance at least $3$. Every code mentioned in this paper is to be considered projective unless otherwise specified. Given a generator matrix $G\in\K^{k\times n}$ of $\Cc$, we can associate $\Cc$ with the projective set $\Pi_G\subseteq \mathbb{P}^{k-1}(\bar \K)$ whose points correspond to the columns of $G$ and where $\bar \K$ denotes the algebraic closure of $\K$. 

Conversely, let $\Pi$ be an ordered set of points in the projective space $\mathbb{P}^{k-1}(\K)$. We can associate a code with $\Pi=\{P_1,\dots,P_n\}$ as follows. For each point $P_i$ we fix a representative in $\K^k$. Given an integer $a\geq 0$, we denote by $\K[x_1,\dots,x_k]_a$ the space of homogeneous polynomials of degree $a$ over $\K$. Consider the evaluation map on the set of representatives of $\Pi$
\begin{equation*}
\begin{split}
    \mathrm{ev}_a:\,&\K[x_1,\dots,x_k]_a\rightarrow \K^n\\
    &f\mapsto (f(P_1),\dots,f(P_n))\,.
\end{split}
\end{equation*}
We denote the image of $\mathrm{ev}_a$ by $C(\Pi)_a$. This is a linear code and it is called the evaluation code of order $a$ associated to $\Pi$. We remark that every linear code can be regarded as an evaluation code of order $1$, in particular $\Cc$ is equivalent to $C(\Pi_G)_1$ with $G$ generator matrix for $\Cc$. 

\begin{definition}
    A code $\Cc\subseteq\K^n$ is called $D$-iso-dual if $\Cc=\Cc^{\perp}D$, where $D$ is a nonsingular diagonal matrix. A code $\Cc$ is called iso-dual if and only if it is monomially equivalent to its dual, i.e., $\Cc$ and $\Cc^{\perp}$ coincide up to coordinate permutations and nonzero coordinate scalings.
\end{definition}
There is a natural correspondence between self-associated sets of points and $D$-iso-dual codes. Indeed, if $\Cc$ is a $D$-iso-dual code with generator matrix $G$, then the set of points $\Pi_G$ is self-associated. Conversely, if $\Pi$ is a self-associated set of points, then a $D$-iso-dual code is obtained directly from the definition. 

Note that every $D$-iso-dual code is iso-dual, but the converse is not always true as one can see in the following example.
\begin{example}\label{ex:isodual}
    Let $\Cc\subseteq\F_7^6$ be the code with generator matrix
\begin{equation*}
    G=(\mathrm{Id}\mid A)=\begin{pmatrix}
        1&0&0&1&1&2\\
        0&1&0&1&3&5\\
        0&0&1&1&5&4
    \end{pmatrix},
\end{equation*}
and let $P$ be the permutation matrix exchanging coordinates $1$ and $4$. A direct computation shows that $G\,(GP)^T=0$, so that $\Cc^\perp=\Cc P$ and $\Cc$ is iso-dual. Moreover, the only solution to the system $GDG^T=0$ is $D=0$. Hence $\Cc$ is not $D$-iso-dual and so $\Pi_G$ is not self-associated.
\end{example}
\begin{remark}\label{remark:algcl}
    If $\K$ is an algebraically closed field, every $D$-iso-dual code is equivalent to a self-dual code. This is not true over an arbitrary field. In fact, if the field is not closed under the square root operation, it may not be possible to decompose $GDG'^T$ as $(G\sqrt{D})(G'\sqrt{D})^T$ in Definition~\ref{def:Gale}. We refer also to~\cite[Lemma 2.1]{betti2025mukai}.
\end{remark}
We denote by $\nu_d$ the Veronese map of degree $d$ that maps a point $P$ in $\mathbb{P}^{k-1}(\K)$ to a point $\nu_d(P)$ in $\mathbb{P}^{\binom{d+k-1}{d}-1}(\K)$, whose coordinates are all the monomials of degree $d$ in the entries of $P$. The $d$-th Veronese code of order $a$ associated with a set $\Pi$ is the evaluation code $C(\nu_d(\Pi))_a$. We recall the following proposition.
\begin{proposition}[{\cite[Proposition 5.66]{tohuaneanu2024commutative}}]
    Let $\Pi$ be a finite set of points and $d,a\geq1$ be integers. Then the $d$-th Veronese code $C(\nu_d(\Pi))_a$ of order $a$ and the evaluation code $C(\Pi)_{ad}$ are equivalent.
\end{proposition}
An arithmetically Gorenstein set of points is a finite set of points in projective space whose homogeneous coordinate ring is Gorenstein, see \cite{brunsherzog93cm}, \cite[Page 215]{Eis2005arit}. The minimum distance of an evaluation code associated to an arithmetically Gorenstein set of points can be bounded from below by the Castelnuovo–Mumford regularity of the coordinate ring as shown in the next proposition.
\begin{proposition}[{\cite[Proposition 3.1]{Tohaneanu2009distance}}]
    If $\Pi\subseteq\mathbb{P}^{k-1}(\K)$ is a nondegenerate arithmetically Gorenstein finite set of points, then for any $1\leq a\leq\reg(\Pi)-1$, one has
    \begin{equation*}
        d(C(\Pi)_a)\geq \reg(\Pi)-a+1\,,
    \end{equation*}
    where $d(C(\Pi)_a)$ is the minimum distance of $C(\Pi)_a$.
\end{proposition}

The following theorem by Eisenbud and Popescu gives us a characterization of an arithmetically Gorenstein set of points.
\begin{theorem}[{\cite[Theorem 7.3]{eisenbud2000projective}}]\label{theorem:eisenbudpopescu}
        If $\Pi\subseteq\mathbb{P}^{k-1}$ is a nondegenerate set of $2k$ points over an algebraically closed field $\K$, then $\Pi$ is arithmetically Gorenstein if and only if $\Pi$ is self-associated and fails by $1$ to impose independent conditions on quadrics.
\end{theorem}

\begin{definition}
    Let $\Pi\subseteq\mathbb{P}^{k-1}$ be a self-associated set of $2k$ points over an algebraically closed field $\K$. If $\Pi$ fails by $m$ to impose independent conditions on the quadrics, we define the Gorenstein defect of $\Pi$ as 
    \begin{equation*}
        \mathrm{gd(\Pi)} = m-1.
    \end{equation*}
    Note that $\mathrm{gd(\Pi)} = 0$ if and only if $\Pi$ is arithmetically Gorenstein, thus the defect is a measure of how far we are from that situation.
\end{definition}
The following question was posed in~\cite[Section 5.5.1]{tohuaneanu2024commutative}.
\begin{question}\label{question}
    If $\Cc\subseteq\K^{2k}$ is a projective self-dual code, does it imply that $\Pi_G$ is arithmetically Gorenstein?
\end{question}
The answer to this question (see Corollary \ref{corollary: mainequivalence}) appears to be folklore among specialists, and a proof can be obtained by combining~\cite[Corollary 2.42]{randriam15course} with Corollary~\ref{corollary:eisenbud} below, see Remark~\ref{remark:randriam1}. As we have mentioned in the introduction, an algebraic proof was also recently given by Kreuzer~\cite[Theorem 4.4]{kreuzer2026codeequivalencepointset}. In this paper, we take a different route and provide a graph-theoretical proof.
We begin by expressing the conditions of Theorem~\ref{theorem:eisenbudpopescu} in terms of the dimension of the second Schur power of a code. The Schur product, or component-wise product, of two vectors $(v_1,\dots,v_n)$ and $(w_1,\dots,w_n)$ in $\K^n$ is defined as
\begin{equation*}
    (v_1,\dots,v_n)\ast(w_1,\dots,w_n)=(v_1  w_1,\dots,v_n w_n)\,.
\end{equation*}
The Schur product $\Cc_1\ast\Cc_2$ of two codes $\Cc_1$ and $\Cc_2$ is the vector space generated by all possible products of a codeword of $\Cc_1$ with a codeword of $\Cc_2$. The $d$-th Schur power $\Cc^{(d)}$ of $\Cc$ is just the Schur product of $d$ copies of $\Cc$. 
\begin{remark}
    Let $\Cc$ be a projective $[n,k]$ linear code with generator matrix $G$. Then, the $d$-th Veronese code $C(\nu_d(\Pi_G))_1$ of order $1$ is equivalent to the $d$-th Schur power $\Cc^{(d)}$.
\end{remark}
Let $\Cc$ be a projective self-dual code of dimension $k$ in $\K^{2k}$ and assume that the set $\Pi_G$ in $\mathbb{P}^{k-1}(\bar\K)$ fails by $m$ to impose independent conditions on the quadrics. This implies that the dimension of the space of quadrics in $\bar\K[x_1,\dots,x_k]_2$ vanishing on $\Pi_G$ is $\binom{k+1}{2}-(2k-m)=\binom{k+1}{2}-(2k-1)+\mathrm{gd}(\Pi_G)$, and so we obtain
\begin{equation}\label{eq:dimdefect}
    \dim(\Cc^{(2)})=2k-\mathrm{gd}(\Pi_G)-1.
\end{equation}
By Theorem~\ref{theorem:eisenbudpopescu}, $\Pi_G$ is arithmetically Gorenstein if and only if $m=1$, equivalently, if and only if $\dim_{\K}(\Cc^{(2)})=2k-1$.
\begin{corollary}\label{corollary:eisenbud}
    Let $\Cc$ be a $k$-dimensional projective self-dual code. Then, $\Pi_G$ is arithmetically Gorenstein if and only if $\dim(\Cc^{(2)})=2k-1$.
\end{corollary}

We are now in position to construct a counterexample to Question~\ref{question}. 
\begin{example}\label{ex: toyexample 1}
Given an $[n_1,k_1]$ code $C_1$ and an $[n_2,k_2]$ code $C_2$, their direct sum is the $[n_1+n_2,k_1+k_2]$ code \[ C_1\oplus C_2 = \{(c_1,c_2):c_1\in C_1,\ c_2\in C_2\}, \] where $(c_1,c_2)$ denotes the concatenation of $c_1$ and $c_2$. Consider the self-dual code $\Cc\subseteq \F_2^8$ defined by the matrix
\begin{equation*}
    G=\begin{pmatrix}
        1&0&0&0&0&1&1&1\\
        0&1&0&0&1&0&1&1\\
        0&0&1&0&1&1&0&1\\
        0&0&0&1&1&1&1&0
    \end{pmatrix}\,.
\end{equation*}
The code $\Cc\oplus\Cc$, generated by the matrix
\begin{equation*}
    G_2=\begin{pmatrix}
        G&0\\
        0&G
    \end{pmatrix}\,
\end{equation*}
is also self-dual, the set $\Pi_{G_2}$ has cardinality $16$, but it fails by $2$ to impose independent conditions on quadrics. In fact, by simple computations one can check that $\dim((\Cc\oplus\Cc)^{(2)})=14$. 
\end{example}

The choice of $\Cc$ in the previous example is not special. Indeed, one can do the same with any other projective self-dual code. In the next section, we give a formal proof of this fact and we will also show that this is essentially the only way to construct self-dual codes whose associated set of points fails by more than $1$ to impose independent conditions on quadrics.

\section{The dimension of the square of a self-dual code and a combinatorial description of the Gorenstein defect}\label{Section 3}
    We recall the definition of indecomposable codes as given in~\cite[Section 1.7]{slepian1960some}.

    \begin{definition}
        A code is said to be decomposable if it is equivalent to the direct sum of two nontrivial codes. A code is called indecomposable if it is not decomposable. %Every code admits an equivalent expression as a sum of indecomposable subcodes. The number of indecomposable subcodes (blocks) in the decomposition of a code $\Cc$ is an invariant of the code, which we denote by $\mathrm{nb}(\Cc)$. 
       Every code is equivalent to a direct sum of indecomposable codes. The number of summands in such a decomposition is an invariant of $\Cc$, called the number of blocks of $\Cc$ and denoted by $\operatorname{nb}(\Cc)$.
        
    \end{definition}
        Given a $k$-dimensional code $\Cc\subseteq\K^n$ the invariant $\mathrm{nb}(\Cc)$ can be found in polynomial time in $n$, by computing the number of connected components of a certain graph. Up to equivalence, let $G=(\mathrm{Id}\,| A)$ be a generator matrix of $\Cc$ in reduced row echelon form. We define the matrix $\bar A=(\bar a_{i,j}) \in\Z^{k\times (n-k)}$ by
        \begin{equation*}
            \bar a_{i,j}=\begin{cases}
                0&\text{if }a_{i,j}=0\\
                1&\text{otherwise.}
            \end{cases}
        \end{equation*}
        The zero-one symmetrization matrix of $A$ is the symmetric matrix
        \begin{equation*}
            \tilde A=\begin{pmatrix}
                0&\bar A\\
                \bar A^t&0
            \end{pmatrix}.\,
        \end{equation*}
        We denote by $\Gamma_{\Cc}$ the bipartite graph with adjacency matrix $\tilde A$. Recall that a symmetric matrix is called irreducible if and only if the associated graph is connected. Although the graph depends on the chosen systematic generator matrix, its number of connected components is independent of this choice and equals $\mathrm{nb}(\Cc)$. Indeed, if $\Cc$ is equivalent to $\Cc_1\oplus\dots\oplus \Cc_{\mathrm{nb}(\Cc)}$, then every systematic generator matrix of $\Cc$ is block diagonal up to row and column permutations, and $\Gamma_{\Cc}$ has at least $\mathrm{nb}(C)$ components. Conversely, permuting rows and columns according to the components of $\Gamma_{\Cc}$ writes $A$, and hence $G$, in block diagonal form with one block per component, so $\Cc$ is a direct sum of $c$ codes and $\mathrm{nb}(\Cc)\ge c$. In particular, $\Cc$ is indecomposable if and only if $\Gamma_{\Cc}$ is connected, that is, if and only if $\tilde A$ is irreducible. The number of connected components can be found in $O(V+E)$ using Tarjan's algorithm, \cite{Tarjan72connected}, where $V$ is the number of vertices and $E$ is the number of edges. In our case $V=n$ and $E\leq (n-k)k< n^2$. Thus, the number of blocks can be computed in $O(n^2)$ time assuming that the systematic matrix is given.
    
    \begin{definition}\label{def: connected sets}
        Let $A$ be an $m\times n$ matrix, with rows $a_1,\ldots,a_m$. We say that a set of pairs $Y=\{\{k_1,\ell_1\},\dots,\{k_t,\ell_t\}\}$ with $1\leq k_i\leq n$, $1\leq \ell_i\leq n$, and $k_i\neq \ell_i$ is connected for $A$ if for each $i$ there exists a row $a_j$ such that $\{k_i,\ell_i\}\subseteq\mathrm{supp}(a_j)$ and for every nonempty $Y_1$ and $Y_2$ such that $Y= Y_1 \cup Y_2$, there exist $\{\hat k,\hat \ell\}\in Y_1$ and $\{\bar k,\bar \ell\}\in Y_2$ such that $\{\hat k,\hat \ell\} \cap \{\bar{k},\bar{\ell}\} \neq \emptyset$.
        % Note that being connected in the sense of the previous definition for $A$ means being connected as vertices in the corresponding bipartite graph associated to the zero-one symmetrization
    \end{definition}
    In other words, $Y$ is connected for $A$ if every pair in $Y$ is contained in the support of a single row of $A$ and the graph whose vertices are the indices appearing in $Y$ and whose edges are the pairs of $Y$ is connected. Indeed, the second condition in the definition says precisely that the edges of this graph cannot be split into two nonempty families with no vertex in common.
   \begin{remark} Let $\Gamma_\Cc$ be the bipartite graph associated with $\tilde A$, whose row vertices and column vertices are, respectively, $\{v_1,\ldots,v_m\}$, $\{w_1,\ldots,w_n\}$. Each pair $\{k,\ell\}$ in a connected set for $A$ corresponds to a path $$ w_kv_jw_\ell $$ in $\Gamma_{\Cc}$ for some $j\in\{1,\ldots,m\}$. Therefore, if $Y$ is connected for $A$, then all column vertices indexed by 
   \begin{equation*}
       \bigcup_{\{k,\ell\}\in Y}\{k,\ell\}
   \end{equation*} belong to the same connected component of $\Gamma_\Cc$. \end{remark}

\begin{example}\label{ex: toyexample 2} Consider the generator matrix $G=(I_4\mid A)$ from Example \ref{ex: toyexample 1}, where 
   \begin{equation*}
       A= \begin{pmatrix} 
       0&1&1&1\\ 
       1&0&1&1\\ 
       1&1&0&1\\ 
       1&1&1&0 
       \end{pmatrix}.
   \end{equation*} Then $$ Y=\bigl\{\{1,3\},\{1,4\},\{3,4\}\bigr\} $$ is a connected set for $A$. Indeed, every pair in $Y$ is contained in $$\operatorname{supp}(a_2)=\{1,3,4\}, $$ and the graph whose edges are the elements of $Y$ is connected. In the bipartite graph $\Gamma_\Cc$ associated with $\tilde A$, see Figure \ref{fig: figexample1}, the vertices $\{5,7,8\}$, corresponding to columns $\{1,3,4\}$ of $A$, are all adjacent to the vertex $2$, which corresponds to the second row of $A$. 

    \begin{figure}[h!]
    \centering
    \begin{tikzpicture}[scale=0.6, every node/.style={circle, draw, minimum size=3mm, inner sep=2pt}]
    % Nodos superiores
    \node (1) at (-4,1) {1};
    \node (2) at (-2,1) {2};
    \node (3) at (0,1) {3};
    \node (4) at (2,1) {4};
    
    % Nodos inferiores
    \node (5) at (-4,-2) {5};
    \node (6) at (-2,-2) {6};
    \node (7) at (0,-2) {7};
    \node (8) at (2,-2) {8};
    
    % Aristas diagonales
    \draw (1) -- (6);
    \draw (1) -- (7);
    \draw (1) -- (8);
    \draw (2) -- (5);
    \draw (2) -- (7);
    \draw (2) -- (8);
    \draw (3) -- (5);
    \draw (3) -- (6);
    \draw (3) -- (8);
    \draw (4) -- (5);
    \draw (4) -- (6);
    \draw (4) -- (7);
    \end{tikzpicture}
    \caption{$\Gamma_{\Cc}$}
    \label{fig: figexample1}
    \end{figure}

    The zero-one symmetrization matrix of $A$ is 
    \begin{equation*}
        \tilde A=\begin{pmatrix}
        0&0&0&0&0&1&1&1\\
        0&0&0&0&1&0&1&1\\
        0&0&0&0&1&1&0&1\\
        0&0&0&0&1&1&1&0 \\
        0&1&1&1&0&0&0&0\\
        1&0&1&1&0&0&0&0\\
        1&1&0&1&0&0&0&0\\
        1&1&1&0&0&0&0&0\\
        \end{pmatrix}\,
    \end{equation*}
   \end{example}

%     \begin{example}\label{ex: toyexample2}
%         Consider the matrix $G$ of Example \ref{ex: toyexample 1}, that defined a self-dual code over $\F_2^8$:

%         \begin{equation*}
%         G=\begin{pmatrix}
%         1&0&0&0&0&1&1&1\\
%         0&1&0&0&1&0&1&1\\
%         0&0&1&0&1&1&0&1\\
%         0&0&0&1&1&1&1&0
%         \end{pmatrix} = (\mathrm{Id} \mid A).
%         \end{equation*}

% The zero-one symmetrization matrix of $A$ is 
% \begin{equation*}
%         \tilde A=\begin{pmatrix}
%         0&0&0&0&0&1&1&1\\
%         0&0&0&0&1&0&1&1\\
%         0&0&0&0&1&1&0&1\\
%         0&0&0&0&1&1&1&0 \\
%         0&1&1&1&0&0&0&0\\
%         1&0&1&1&0&0&0&0\\
%         1&1&0&1&0&0&0&0\\
%         1&1&1&0&0&0&0&0\\
%         \end{pmatrix}\,
% \end{equation*}

%     Denote by $\{a_1, \ldots,a_8\}$ the rows of $\tilde A$. The set $Y = \{\{5,7\}, \{5,8\}\,\{7,8\}\}$ is a connected set for $\tilde A$. Observe that $\{5,7\}, \{5,8\},\{7,8\} \subseteq \mathrm{supp}(a_2)$. In Figure \ref{fig: fig1}, the vertices $\{5,7,8\}$ are connected to the vertex $\{2\}$.

%     \end{example}
    The next lemma says that $\Gamma_{\Cc}$ is connected exactly when all column indices of $A$ can be linked to one another through a chain of pairs, each of which lies in the support of a common row. This lets us replace the connectivity of the bipartite graph $\Gamma_{\Cc}$ by a condition involving only the supports of the rows of $A$.
    \begin{lemma}\label{lemma:irrimpliesconnected}
        Let $n\geq2$ and let $A$ be an $m\times n$ matrix with no zero rows. The zero-one symmetrization matrix $\tilde A$ is irreducible if and only if there exists a connected set $Y=\{\{k_1,\ell_1\},\dots,\{k_t,\ell_t\}\}$ for $A$ such that $[n]=\{k_1,\ell_1\}\cup\dots\cup\{k_t,\ell_t\}$.
    \end{lemma}
    %\begin{proof}
    %    By definition, $\tilde A$ is irreducible if and only if the associated graph $\Gamma$ is connected. The graph $\Gamma$ is bipartite by construction. Let $V=\{v_1,\dots,v_m\}$ be the set of vertices whose indices correspond to the rows of $A$ and $W=\{w_1,\dots,w_n\}$ the set of vertices whose indices correspond to the columns of $A$. Every edge of $\Gamma$ goes from $V$ to $W$. Fix a spanning walk $T=t_1,\dots,t_s$ in $\Gamma$, i.e., a walk that passes through every vertex in $\Gamma$. This walk exists since the graph is assumed to be connected. Consider the set
    %   \begin{equation*}
    %        Y=\{\{k,l\}:\,\text{if there exists $i$ such that }t_{i}=w_k\text{ and }t_{i+2}=w_\ell\}\,.
    %   \end{equation*}
    %    Since the walk $T$ is a spanning walk clearly we have $[n]=\bigcup_Y\{k,\ell\}$. Moreover, for every $\{k,l\}\in Y$ there exists $i$ such that $t_{i}=w_k\text{ and }t_{i+2}=w_l$, and so $\{k,l\}$ is contained in the support of the row that corresponds to $t_{i+1}$. Finally, suppose that $Y=Y_1\cup Y_2$, following the walk $T$, there is a moment in which we pass from $Y_1$ to $Y_2$ (or vice versa), which means that $\{t_{i},t_{i+2}\}$ is in $Y_1$ and $\{t_{i+2},t_{i+4}\}$ is in $Y_2$. Therefore, we obtain that $Y$ is connected.

    %    Conversely, suppose that there exists a connected set $Y$. Then, we can construct using $Y$ a connected bipartite graph that is a subgraph of the graph $\Gamma$ with adjacency matrix equal to $\tilde A$. This implies that $\Gamma$ is also connected and so we conclude that $\tilde A$ is irreducible.
    %\end{proof}
    \begin{proof}
        Let $\Gamma$ be the bipartite graph whose adjacency matrix is $\tilde A$, with row vertices $v_1,\ldots,v_m$ and column vertices $w_1,\ldots,w_n$. Suppose first that $\tilde A$ is irreducible, that is, $\Gamma$ is connected. Choose a column vertex $w_{k_0}$ and set $S_0=\{k_0\}$. Assume that a connected set of column vertices $S_{r-1}\subsetneq [n]$ has already been constructed. Since $\Gamma$ is connected, there exist $k\in S_{r-1}$, $\ell \in [n]\setminus S_{r-1}$, and a row vertex $v_j$ such that 
        $$ \{k,\ell\}\subseteq\operatorname{supp}(a_j).$$
        Equivalently,
        $$ w_kv_jw_\ell $$ is a path in $\Gamma$. Add the pair $\{k,\ell\}$ to $Y$ and set $S_r = S_{r-1} \cup \{\ell\}$. After $n-1$ steps, $S_{n-1} = [n]$. Moreover, each pair added to $Y$ joins a new index to an index already in the set. Therefore, the set $Y$ is a connected set for $A$ and 
        \begin{equation*}
             [n]=\bigcup_{\{k,\ell\}\in Y}\{k,\ell\}.
        \end{equation*}
        Conversely, suppose that such a connected set $Y$ exists. For every $\{k,\ell\}\in Y$, there is a row $a_j$ such that
        $$\{k,\ell\}\subseteq\operatorname{supp}(a_j), $$ and hence $\Gamma$ contains the path $$ w_kv_jw_\ell.$$ The pairs of $Y$ form a connected graph and cover $[n]$, so all column vertices lie in the same connected component. Since every row vertex is adjacent to a column vertex, $\Gamma$ is connected and thus, $\tilde A$ is irreducible.
    \end{proof}

    \begin{example} 
    Let $A$ be the matrix from Example \ref{ex: toyexample 2}, and let $\Gamma_\Cc$ be the bipartite graph associated with its zero-one symmetrization $\tilde A$, as shown in Figure \ref{fig: figexample1}. The vertices $\{1,2,3,4\}$ correspond to the rows of $A$, whereas the vertices $\{5,6,7,8\}$ correspond, respectively, to columns $\{1,2,3,4\}$ of $A$. We illustrate the construction in the proof of Lemma \ref{lemma:irrimpliesconnected}. Start with the third column of $A$ and set
    $$S_0=\{3\}, \qquad Y_0=\emptyset.$$
    In $\Gamma_\Cc$, the column vertices corresponding to columns $3$ and $1$ are joined through the second row vertex: \[ 7-2-5. \] Equivalently, \[ \{3,1\}\subseteq\operatorname{supp}(a_2). \] We therefore set \[ S_1=\{3,1\}, \qquad Y_1=\bigl\{\{3,1\}\bigr\}. \] Next, the column vertices corresponding to columns $1$ and $2$ are joined through the fourth row vertex: \[ 5-4-6. \] Thus, \[ \{1,2\}\subseteq\operatorname{supp}(a_4), \] and we set \[ S_2=\{3,1,2\}, \qquad Y_2=\bigl\{\{3,1\},\{1,2\}\bigr\}. \] Finally, the column vertices corresponding to columns $2$ and $4$ are joined through the third row vertex: \[ 6-3-8. \] Hence, \[ \{2,4\}\subseteq\operatorname{supp}(a_3), \] and we set 
    $$ S_3=\{3,1,2,4\}=[4]$$ and
    $$Y_3= \bigl\{\{3,1\},\{1,2\},\{2,4\}\bigr\}.$$
    Therefore,
    $$Y:=Y_3= \bigl\{\{3,1\},\{1,2\},\{2,4\}\bigr\}$$ is connected for $A$ and satisfies 
    \begin{equation*}
        [4]=\bigcup_{\{k,\ell\}\in Y}\{k,\ell\}. 
    \end{equation*}
    \end{example}

    \begin{lemma}\label{lemma:generatingset}
        Let $A$ be an $m\times n$ matrix, $Y=\{\{k_1,\ell_1\},\dots,\{k_t,\ell_t\}\}$ a connected set for $A$. For each $i \in [t]$, let $c_i=e_{k_i}-e_{\ell_i}\in\K^n$. If $[n]=\{k_1,\ell_1\}\cup\dots\cup\{k_t,\ell_t\}$, then
        \begin{equation*}
            \ker(x_1+\dots+x_n)=\langle c_{1},\dots,c_{t}\rangle_{\K}\,.
        \end{equation*}
    \end{lemma}

    \begin{proof} 
    Since every $c_i$ belongs to $\ker(x_1+\cdots+x_n)$, we have that 
    \begin{equation*}
        \langle c_1,\ldots,c_t\rangle_K \subseteq \ker(x_1+\cdots+x_n).
    \end{equation*}  Regarding the pairs in $Y$ as edges, they form a connected graph on $[n]$. Choose a spanning tree of this graph and fix the vertex $1$. For any $j\in[n]$, that spanning tree contains a unique path joining the vertex $1$ and $j$, 
    \begin{equation*}
        1=i_0,i_1,\ldots,i_r=j.
    \end{equation*}Taking the appropriate signed sum of the vectors corresponding to the edges of the unique path from $1$ to $j$ gives
    \begin{equation*}
        e_1-e_j\in\langle c_1,\ldots,c_t\rangle_{\K}. 
    \end{equation*} Since \[ e_1-e_2,\ldots,e_1-e_n \] form a basis of $\ker(x_1+\cdots+x_n)$, the reverse inclusion follows. 
    \end{proof}

    \begin{example} 
    Consider the matrix $A$ from Example \ref{ex: toyexample 2}: 
    \begin{equation*}
         A= 
         \begin{pmatrix} 
         0&1&1&1\\ 
         1&0&1&1\\ 
         1&1&0&1\\ 
         1&1&1&0 
         \end{pmatrix}.
    \end{equation*}
     The set \[ Y=\bigl\{\{4,3\},\{3,1\},\{1,2\}\bigr\} \] is a connected set for $A$ and its pairs cover $[4]$. The corresponding vectors in Lemma \ref{lemma:generatingset} are $$ c_1=e_4-e_3=(0,0,1,1), $$ $$ c_2=e_3-e_1=(1,0,1,0), $$ and $$ c_3=e_1-e_2=(1,1,0,0). $$ Therefore, $$ \ker(x_1+x_2+x_3+x_4) = \langle c_1,c_2,c_3\rangle_{\mathbb F_2}. $$
    \end{example}

    \begin{proposition}\label{proposition:2k-1}
        Let $\Cc\subseteq\K^{2k}$ be an indecomposable self-dual code. Then,
        $$\dim(\Cc^{(2)})=2k-1\,.$$
    \end{proposition}
    \begin{proof}
        Up to a permutation of coordinates, let $(\mathrm{Id} \mid A)$ with $A$ such that $AA^T=-\mathrm{Id}$ be a generator matrix of $\Cc$, and let $a_1,\dots,a_k$ be the rows of $A$. Then, clearly $\dim(\Cc^{(2)})=2k-1$ if and only if $\dim(\mathcal{A})=k-1$, where
        $$\mathcal{A}=\langle\{a_i\ast a_j:i\neq j\}\rangle_{\K}\,.$$
        Since the code is self orthogonal, we have that
        $$\mathcal{A}\subseteq \ker(x_{1}+\dots+x_k)\,.$$
        Since the dimension of $\ker(x_{1}+\dots+x_k)$ is $k-1$, to conclude we have to prove that the previous inclusion is in fact an equality. By Lemma~\ref{lemma:irrimpliesconnected}, there exists a connected set of pairs $Y=\{\{s_1,\ell_1\},\dots,\{s_t,\ell_t\}\}$ such that $[k]=\{s_1,\ell_1\}\cup\dots\cup\{s_t,\ell_t\}$. Given two entries $s,\ell$ up to a nonzero constant there exists only one element in $\ker(x_{1}+\dots+x_k)$ with support equal to $\{s,\ell\}$. We call this element $c_{s,\ell}=e_s-e_\ell$. By Lemma~\ref{lemma:generatingset} the vector space $\ker(x_{1}+\dots+x_k)$ is generated by the set of all $c_{s,\ell}$ with $\{s,\ell\}\in Y$. By definition of connected set, there exists a row $a_i$ such that $\{s,\ell\}\subseteq \mathrm{supp}(a_i)$. Without loss of generality let $i=1$. The vector space $a_1\ast \langle a_2,\dots, a_k\rangle_{\K}$ is contained in $\mathcal{A}$ and it has dimension at least $\lvert \mathrm{supp}(a_1)\rvert-1$. Indeed, consider the map $\pi:\langle a_2,\dots,a_k\rangle_{\K}\rightarrow \K^{\lvert\mathrm{supp}(a_1)\rvert}$ given by multiplying by $a_1$ and then by projecting on its support. We have $\dim(\ker \pi)\leq k-\lvert\mathrm{supp}(a_1)\rvert$, and so
        $$\dim(a_1\ast \langle a_2,\dots, a_k\rangle_{\K})+k-\lvert\mathrm{supp}(a_1)\rvert\geq\dim(a_1\ast \langle a_2,\dots, a_k\rangle_{\K})+\dim(\ker\pi)=k-1,$$
        where in the last equality we used the fact that $\dim \langle a_2,\dots,a_k\rangle_{\K}=k-1$ since $A$ is invertible.
        Let $\ker(x_{1}+\dots+x_k)_{a_1}$ be the subspace of $\ker(x_{1}+\dots+x_k)$ supported on the support of $a_1$. Then, 
        $$\dim_{\K}(\ker(x_{1}+\dots+x_k)_{a_1})=\lvert \mathrm{supp}(a_1)\rvert-1\text{ and }a_1\ast \langle a_2,\dots, a_k\rangle_{\K}\subseteq \ker(x_{1}+\dots+x_k)_{a_1}\,,$$
        and so, by a dimension argument we conclude that 
        $$a_1\ast \langle a_2,\dots, a_k\rangle_{\K}= \ker(x_{1}+\dots+x_k)_{a_1}\,.$$
        Finally, 
        $$c_{s,\ell}\in \ker(x_{1}+\dots+x_k)_{a_1}=a_1\ast \langle a_2,\dots, a_k\rangle_{\K}\subseteq \mathcal{A}\,.$$
        Since this is true for every pair $\{s,\ell\}\in Y$, we obtain $\mathcal{A}=\ker(x_{1}+\dots+x_k)$, and this concludes the proof.
    \end{proof}
    \begin{example}\label{ex: toyexample 3}
        Consider the code $\Cc$ defined in Example \ref{ex: toyexample 1} by the generator matrix $G$:\begin{equation*}
        G=\begin{pmatrix}
        1&0&0&0&0&1&1&1\\
        0&1&0&0&1&0&1&1\\
        0&0&1&0&1&1&0&1\\
        0&0&0&1&1&1&1&0
        \end{pmatrix} = (\mathrm{Id} \mid A)\,
        \end{equation*}
        Denote the rows of $G$ by $\{g_1,g_2,g_3,g_4\}$ and let $e_i$ be the $i$-th standard vector. We have that $\Cc^{(2)} = \langle g_i\ast g_j \, \colon 1\leq i \leq j \leq 4\rangle$, namely:
    \begin{equation*}
        \Cc^{(2)} = \langle g_1,g_2,g_3,g_4, e_7+e_8, e_6+e_8,e_6+e_7,e_5+e_8,e_5+e_7,e_5+e_6\rangle
    \end{equation*}
    A direct computation shows that $\dim(\Cc^{(2)}) = 7 = 2k-1.$
    \end{example}
    The next result generalizes Proposition \ref{proposition:2k-1} stating that the dimension of the Schur product of a self-dual code $\Cc$ can be computed in terms of the number of blocks of $\Cc$.
    \begin{theorem}\label{maintheorem}
        Let $\Cc\subseteq\K^{2k}$ be a self-dual code. Then,
        $$\dim(\Cc^{(2)})=2k-\nbl(\Cc)\,.$$
    \end{theorem}
    \begin{proof}
        Let $\Cc$ be a $k$-dimensional self-dual code and let $\Cc_1,\dots,\Cc_{\nbl(\Cc)}$ be indecomposable codes such that after a permutation of coordinates $\Cc$ is equal to $\Cc_1\oplus\dots\oplus\Cc_{\nbl(\Cc)}$. Since
        \begin{equation*}                \left(\Cc_1\oplus\dots\oplus\Cc_{\nbl(\Cc)}\right)^{(2)}=\Cc_1^{(2)}\oplus\dots\oplus\Cc_{\nbl(\Cc)}^{(2)}\,,
        \end{equation*}
        we obtain
        \begin{equation*}
            \dim\left(\Cc^{(2)}\right)=\dim\left(\Cc_1^{(2)}\right)+\dots+\dim\left(\Cc_{\nbl(\Cc)}^{(2)}\right)\,.
        \end{equation*}
        By applying Proposition~\ref{proposition:2k-1} for each $\Cc_{i}$, since the indecomposable components of a self-dual code are self-dual, we obtain
        \begin{equation*}
            \dim\left(\Cc^{(2)}\right)=\sum_{i=1}^{\nbl(\Cc)}(2k_i-1)=2k-\nbl(\Cc)\,,
        \end{equation*}
        where $k_i$ is the dimension of $\Cc_i$ for all $i\in[\nbl(\Cc)]$.
    \end{proof}
    \begin{remark}
        Note that the previous theorem remains valid for $D$-iso-dual codes. In fact, both the dimension of $\Cc^{(2)}$ and the number of blocks $\nbl(\Cc)$ are invariant under field extension, and we already observed that over an algebraically closed field every $D$-iso-dual code is equivalent to a self-dual code. However, it fails for iso-dual codes in general. Indeed, let $\Cc$ be the code of Example~\ref{ex:isodual}. A direct computation shows that $\dim \Cc^{(2)}=6$, whereas $n-\nbl(\Cc)=5$ since $\Cc$ is indecomposable. This also shows that the conclusion of \cite[Theorem 4.4]{kreuzer2026codeequivalencepointset} does not hold for iso-dual codes in general.
    \end{remark}
    \begin{remark}\label{remark:randriam1}
        The previous theorem can be derived as a special case of \cite[Corollary 2.42]{randriam15course}. There, the author defines the extended stabilizing algebra of a linear code as
        \begin{equation*}
            \widehat{\mathcal{A}}=\{a \in \K^n: a \ast C\subseteq C \}
        \end{equation*} and it was proven that \cite[Proposition 2.11, Corollary 2.42]{randriam15course}:
        \begin{equation*}
            \dim(\widehat{\mathcal{A}}) = |\mathcal{P}(C)|,\quad  \widehat{\mathcal{A}} = (\Cc \ast \Cc^\bot)^\bot,
        \end{equation*} where $|\mathcal{P(C)}|$ is the number of components of the code. However, the approach in the present work provides an effective method of computation by exploiting combinatorial methods.
    \end{remark}
   Theorem \ref{maintheorem} gives an explicit formula of the Gorenstein defect in terms of the blocks of the code.
   \begin{corollary}\label{corollary: gd formula}
       Let $\Cc \subseteq \K^{2k}$ be a projective self-dual code with generator matrix $G$. Then,
        \begin{equation*}
            \mathrm{gd}(\Pi_G) = 2k - \dim(\Cc^{(2)})-1 = \nbl(\Cc)-1.
        \end{equation*}
   \end{corollary}

   \begin{proof}
       Follows from Equation~\eqref{eq:dimdefect} and Theorem~\ref{maintheorem}.
   \end{proof}

       \begin{corollary}\label{corollary: mainequivalence}
        Let $\Cc\subseteq\K^{2k}$ be a projective self-dual code. Then, $\Pi_G$ is arithmetically Gorenstein if and only if $\Cc$ is indecomposable.
    \end{corollary}
    \begin{proof}
        Follows immediately from Theorem~\ref{maintheorem} and Corollary~\ref{corollary:eisenbud}.
    \end{proof}
Corollary \ref{corollary: mainequivalence} implies directly that the column points of the code in Example \ref{ex: toyexample 1} are arithmetically Gorenstein, a fact that we already verified by computing the dimension of the Schur square in Example \ref{ex: toyexample 3}. Moreover, we can compute the vanishing ideal of the points to confirm this.
    
    \begin{example}
        Consider the self-dual code $\Cc \subseteq \F_2^8$ of Example \ref{ex: toyexample 1}, whose generator matrix is: 
        \begin{equation*}
        G=\begin{pmatrix}
        1&0&0&0&0&1&1&1\\
        0&1&0&0&1&0&1&1\\
        0&0&1&0&1&1&0&1\\
        0&0&0&1&1&1&1&0
        \end{pmatrix}
        \end{equation*}
        We already know that $\Cc$ is indecomposable, so $\Pi_G$ must be Gorenstein. The ideal of the column points $\Pi_{G}$ is:
        \begin{equation*}
            I = \langle x_3(x_0+x_1+x_2),~ x_2(x_0+x_1+x_3),~x_1(x_0+x_2+x_3)\rangle,
        \end{equation*} which is a Gorenstein ideal since it is a complete intersection of three quadrics. Indeed, the space of the quadrics in $R=\bar{\F}_2[x_0,\ldots,x_3]$ has dimension $\binom{5}{2}=10$, the Hilbert function of the coordinate ring $R/I$ at degree $2$ gives us the number of independent conditions that the points impose on quadrics, which in this case is $\mathrm{HF}_{R/I}(2) = 7 = (2k-1)-\mathrm{gd}(\Pi_G)$, so the Gorenstein defect is $0$, as we already knew because the points are arithmetically Gorenstein since the matrix $\tilde A$ is irreducible. Then, $I_2$ has dimension $\binom{5}{2}-\mathrm{HF}_{R/I}(2) = 10-7 = 3$.
    \end{example}

We now consider Example \ref{ex: toyexample 1}, our first counterexample to Toh{\u{a}}neanu's question. Corollary \ref{corollary: gd formula} states that the column points of the code in Example \ref{ex: toyexample 1} with generator matrix $G_2$ have Gorenstein defect $1$. We compute the ideal of these column points and verify that Corollary \ref{corollary: gd formula} holds.

    \begin{example}
Consider the code $\Dd = \Cc \oplus \Cc$ from Example \ref{ex: toyexample 1}, with generator matrix $G_2$. The Gorenstein defect of the set of column points is $1$, since the code has two blocks: $\mathrm{nb}(\Dd)-1 = 2-1 =1$. Moreover, notice that the associated graph $\Gamma_{\Dd}$ of the zero-one symmetrization matrix has two connected components, as shown in Figure \ref{fig: fig4}.

\begin{figure}[h!]
    \centering
    \begin{tikzpicture}[scale=0.6, every node/.style={circle, draw, minimum size=5mm, inner sep=0pt}]
    % Nodos superiores
    \node (1) at (-6,1) {1};
    \node (2) at (-4,1) {2};
    \node (3) at (-2,1) {3};
    \node (4) at (0,1) {4};
    \node (5) at (2,1) {5};
    \node (6) at (4,1) {6};
    \node (7) at (6,1) {7};
    \node (8) at (8,1) {8};
    
    % Nodos inferiores
    \node (9) at (-6,-2) {9};
    \node (10) at (-4,-2) {10};
    \node (11) at (-2,-2) {11};
    \node (12) at (0,-2) {12};
    \node (13) at (2,-2) {13};
    \node (14) at (4,-2) {14};
    \node (15) at (6,-2) {15};
    \node (16) at (8,-2) {16};
    
    % % Aristas del 1
    \draw (1) -- (10);
    \draw (1) -- (12);
    \draw (1) -- (11);
    % %Aristas del 2
    \draw (2) -- (9);
    \draw (2) -- (11);
    \draw (2) -- (12);
    % %Aristas del 3
    \draw (3) -- (9);
    \draw (3) -- (10);
    \draw (3) -- (12);
    % %Aristas del 4
    \draw (4) -- (9);
    \draw (4) -- (10);
    \draw (4) -- (11);
    % %Aristas del 5
    \draw (5) -- (14);
    \draw (5) -- (15);
    \draw (5) -- (16);
    % %Aristas del 6
    \draw (6) -- (15);
    \draw (6) -- (16);
    \draw (6) -- (13);
    % %Aristas del 7
    \draw (7) -- (14);
    \draw (7) -- (16);
    \draw (7) -- (13);
    % %Aristas del 8
    \draw (8) -- (15);
    \draw (8) -- (14);
    \draw (8) -- (13);
    \end{tikzpicture}
    \caption{$\Gamma_{\Dd}$}
    \label{fig: fig4}
    \end{figure}

The ideal $J$ of the column points of $\Dd$ is generated by the following $22$ quadrics:
 
\begin{align*}
    J = \langle & x_4x_7+x_5x_7+x_6x_7,~x_3x_7,~ x_2x_7,~ x_1x_7,~ x_0x_7,~x_4x_6+x_5x_6+x_6x_7,~x_3x_6,\\ &x_2x_6,~x_1x_6,~x_0x_6,~x_3x_5,~x_2x_5,~x_1x_5,~x_0x_5,~x_3x_4,~x_2x_4,~x_1x_4,~x_0x_4,\\ & x_0x_3+x_1x_3+x_2x_3,~x_0x_2+x_1x_2+x_2x_3,~x_0x_1+x_1x_2+x_1x_3,x_4x_5+x_5x_6+x_5x_7\rangle
\end{align*}

The number of independent conditions that the column points impose on quadrics is equal to the Hilbert function of the coordinate ring $R/J$ evaluated at $2$, which in this case is $\mathrm{HF}_{R/J}(2) = 14$. Thus, the $16$ column points fail by $2$ to impose independent conditions in the space of quadrics $\bar{\F}_2[x_0, \ldots,x_7]_2$, which has dimension $\binom{9}{2} = 36$. Hence, the Gorenstein defect is 1, as expected.
\end{example}

We remark that although our approach comes from coding theory, the results actually work over any field: the graph we construct and all the techniques are field independent. In particular, the combinatorial criterion via Schur square dimension and connected components gives a direct way to check when  self-associated points over an algebraically closed field are arithmetically Gorenstein (failing by one on quadrics).

\begin{corollary}\label{corollary:alg closed}
    Let $\K$ be an algebraically closed field. Let $\Pi\subset\PP^{k-1}(\K)$ be a nondegenerate self-associated ordered set of $2k$ distinct points. Then $\Pi$ is arithmetically Gorenstein if and only if $\Cc(\Pi)_1$ is indecomposable.
\end{corollary}
\begin{proof}
 It follows from Remark~\ref{remark:algcl} and Corollary~\ref{corollary: mainequivalence}.
\end{proof}
 
\section{Number of Gorenstein self-dual codes}\label{Section 4}

It is shown in Section \ref{Section 3} that self-dual codes whose column points define an arithmetically Gorenstein set are exactly those that are indecomposable, which answers Question \ref{question}. The next question is whether there exist self-dual codes whose columns points define an arithmetically Gorenstein set for any given length $n$ over a fixed finite field $\F_q$, or equivalently by Corollary \ref{corollary: mainequivalence}, whether there are indecomposable self-dual codes for any length $n$.

It was shown in \cite[Th. 4.1.3]{Haluk1989phd} that, for a given length $n$ large enough, almost all binary self-dual codes are indecomposable, or equivalently, have Gorenstein defect equal to $0$. Namely, denote by $\Gg_{n,2}$ the number of self-dual codes of length $n$, ($n$ must be even) over $\F_2$ and $\Cc_{n,2}$ the number of indecomposable self-dual codes of length $n$ over $\F_2$. Then \cite[Th. 4.1.3]{Haluk1989phd} states that
\begin{equation*}
    \lim_{n\to\infty} \frac{\Cc_{n,2}}{\Gg_{n,2}}=1
\end{equation*}

Nevertheless, that result can be extended in a simple way to any $q$. Denote by $\Gg_{n,q}$ the number of self-dual codes of length $n$, ($n$ must be even) over $\F_q$ and $\Cc_{n,q}$ the number of indecomposable self-dual codes of length $n$ over $\F_q$. We set $\Gg_{0,q}=1$. Then, number $\Gg_{n,q}$ can be computed with the following formula:
\begin{equation*}
\mathcal{G}_{n,q}=\varepsilon_{n,q}\prod_{i=1}^{n/2-1}\bigl(q^{i}+1\bigr),
\qquad
\varepsilon_{n,q}=
\begin{cases}
1 & \text{if } q \text{ is even},\\
2 & \text{if } q \text{ is odd and } (-1)^{n/2} \text{ is a square in } \mathbb{F}_q,\\
0 & \text{otherwise},
\end{cases}
\end{equation*}
this follows from the results in \cite{pless1965number} together with the observation that $\sum x_i^2$ admits an $n/2$-dimensional isotropic subspace iff $(-1)^{n/2}$ is a square.
The proof in \cite[Theorem 4.1.3]{Haluk1989phd} uses the following recursive formula, that is valid for any $q$, 
\begin{equation*}\label{eq: counting}
        n\Gg_{n,q} = \sum_{t=1}^{n/2} \binom{n}{2t}2t\,\Cc_{2t,q}\,\Gg_{n-2t,q},
\end{equation*}
The computations for the general case are the same for any $q$ even and just differ by a factor of $2$ if $q$ is odd, which does not affect the asymptotic behaviour. Thus, the result can be generalized to the following theorem.
\begin{theorem}\label{th: asymptotic}
    Let $n$ run over the even integers if $q$ is even or $q\equiv 1 \pmod 4$ and over the multiples of $4$ if $q\equiv 3 \pmod 4$. Then
    \begin{equation*}
    \lim_{n\to\infty} \frac{\Cc_{n,q}}{\Gg_{n,q}}=1.
    \end{equation*}
    In other words, asymptotically almost all self-dual codes over $\F_q$ are indecomposable.
\end{theorem}

\begin{proof}
    We have that $\Gg_{n,q}>0$ for the values of $n$ considered in the statement of the theorem. By counting the pairs $(\mathcal{C},j)$ where $\mathcal{C} \subset \F_q^n$ is a self-dual code and $j \in [n]$ we have that: 
    \begin{equation}\label{eq: recursion}
        n\mathcal{G}_{n,q} = \sum_{t=1}^{n/2} \binom{n}{2t}2t\mathcal{C}_{2t,q}\mathcal{G}_{n-2t,q}.
    \end{equation}
    Indeed, we have $n$ different coordinates and $\Gg_{n,q}$ self-dual codes of length $n$ so the left side is clear. On the other hand, given a pair $(\Cc,j)$, let $B$ the block containing the coordinate $j$, with $|B| = 2t$. For a given $t$, there are $\binom{n}{2t}$ choices of the set $B$ and $2t$ choices of $j \in B$ so there are $\binom{n}{2t}2t$ choices for the pair $(B,j)$. Hence, we have $\binom{n}{2t}2t\Cc_{2t,q}$ choices for the component that contains the coordinate $j$. The remaining $n-2t$ different coordinates outside the component determine a self-dual code of length $n-2t$ so the sum $\sum_{t=1}^{n/2} \binom{n}{2t}2t\mathcal{C}_{2t,q}\mathcal{G}_{n-2t,q}$ also counts the pairs $(\Cc,j)$.

    Dividing both sides of \eqref{eq: recursion} by $n\Gg_{n,q}$ we have that 
    \begin{align*}
        1 = \frac{\Cc_{n,q}}{\Gg_{n,q}}+F_n, \text{ where } F_n=\sum_{t=1}^{n/2-1} \binom{n-1}{2t-1}\frac{\Cc_{2t,q}\Gg_{n-2t,q}}{\Gg_{n,q}} \leq \sum_{t=1}^{n/2-1} \binom{n-1}{2t-1}\frac{\Gg_{2t,q}\Gg_{n-2t,q}}{\Gg_{n,q}}.
    \end{align*}
    Put $R_t = \frac{\Gg_{2t,q}\Gg_{n-2t,q}}{\Gg_{n,q}}$. Since $R_t = R_{n/2-t}$ and $\binom{n-1}{2t-1}+\binom{n-1}{2(n/2-t)-1} = \binom{n}{2t}$, pairing the terms $t$ and $n/2-t$ gives
    \begin{equation*}
        F_n \leq \sum_{t=1}^{\lfloor n/4 \rfloor}\binom{n}{2t}R_t.
    \end{equation*}
    We now bound $R_t$ for $1\leq t\leq n/4$. The factor $\frac{\varepsilon_{2t,q}\varepsilon_{n-2t,q}}{\varepsilon_{n,q}} \leq 2$. Indeed, it equals 1 for $q$ even and it is in $\{0,2\}$ for $q$ odd. Therefore, using that $q^i \leq q^i+1 \leq q^{i+1}$,
    \begin{align*}
        R_t \leq 2 \frac{\Pi_{i=1}^{t-1}q^{i+1}}{\Pi_{i=n/2-t}^{n/2-1}q^i} = 2q^{t-1-t(n/2-t)}\leq q^{-t(n/2-t-1)}\leq q^{-t(n/4-1)}\leq 2^{-t(n/4-1)}, 
    \end{align*}where we have used $q \geq 2$ and $n/2-t\geq n/4$. The case with $q=2$ is the one considered in \cite{Haluk1989phd}, and the factor of $2$ coming from $\varepsilon$ for $q$ odd is absorbed by $q^{-1}$. Since $\binom{n}{2t}\leq n^{2t}$, we get
    \begin{align*}
        0 \leq F_n \leq \sum_{t=1}^{\lfloor n/4 \rfloor}\Big(\frac{n^2}{2^{n/4-1}}\Big)^t
    \end{align*}
    Let $x_n = \Big(\frac{n^2}{2^{n/4-1}}\Big)$.Since $\lim_{n \rightarrow \infty}x_n = 0$, we can choose $n$ large such that $x_n < 1$. Then,
    \begin{align*}
        F_n \leq \sum_{t=1}^{\lfloor n/4 \rfloor}\Big(\frac{n^2}{2^{n/4-1}}\Big)^t \leq \sum_{t=1}^{\infty}\Big(\frac{n^2}{2^{n/4-1}}\Big)^t = \frac{x_n}{1-x_n}.
    \end{align*}
    Since $x_n \rightarrow0$, we obtain $F_n \rightarrow0$, hence $\frac{\Cc_{n,q}}{\Gg_{n,q}}\rightarrow1$.
\end{proof}

\begin{remark}
  Besides the asymptotic behavior described in Theorem \ref{th: asymptotic}, the number of indecomposable self-dual codes of length $n$ can be explicitly computed by the following formula \cite{BOGART1977333}. Suppose we have a partition of $n$ that has $p_i$ parts of size $i$, so the number of parts is $p = \sum_{i=1}^{n}p_i$. Then, the number of indecomposable self-dual codes of length $n$ is given by 
    \begin{equation*}
        \Cc_{n,q} = n! \sum_{\text{partitions of $n$}} (-1)^{(\sum_ip_i)-1}[(\sum_ip_i)-1]! \prod_{j=1}^{n}\frac{\Gg_{j,q}^{p_j}}{(j!)^{p_j}(p_j)!}
    \end{equation*}  
    The previous formula shows that for $n \in \{4,6,10\}$, there are no indecomposable binary self-dual codes of length $n$. 
Nevertheless, if $n=2k$ is large enough and $\Cc \subseteq \K^{2k}$ is a self-dual code, then almost surely $\Pi_G$ is arithmetically Gorenstein.
\end{remark}

Theorem \ref{th: asymptotic} studies the asymptotics of the indecomposable codes among the self-dual codes. However, the results presented in Section \ref{Section 3} work for projective self-dual codes, that is, codes with no proportional columns. The rest of this Section is devoted to study the behavior of projective indecomposable self-dual codes among the self-dual codes.\\

Let $G$ be the generator matrix of a self-dual code $\Cc \subset \F_q^n$ and let $g_1,\ldots,g_n$ be its columns. We have that $g_j = \lambda g_i$ if and only if $\lambda e_i-e_j \in \Cc^{\perp}  =\Cc$. Hence, a code is projective exactly when $\Cc^{\perp}$ has no nonzero word of weight at most 2. Note that for a self-dual code it is impossible to have a weight 1 word so a self-dual code $\Cc$ is projective if and only if $d(\Cc) \geq 3$.

Recall that every self-dual code $\Cc$ decomposes uniquely as 
\begin{equation}\label{eq: decomposition in projectives}
    \Cc = \Dd_1 \oplus \cdots \oplus \Dd_t \oplus \Pp,
\end{equation}where $\Dd_1, \ldots \Dd_t$ are self dual codes of length $2$ and $\Pp$ is a projective self-dual code. Indeed, let $\Cc$ be a self-dual code over $\F_q^n$ and let $x$ be a codeword in $\Cc$ with $\mathrm{supp}(x) = \{i,j\}$, for any $c \in \Cc$, $\langle x,c \rangle = 0$, so the coordinates $c_i,c_j$ lie in the line $\{(u,v) \colon x_iu+x_jv=0\}$, hence $(c_i,c_j) = \lambda(x_i,x_j)$ for some $\lambda \in \F_q$ and $c-\lambda x$ has support off $\{i,j\}$. Hence, $\Cc = \langle x\rangle\oplus \Cc'$, with $\Cc'$ self-dual of length $n-2$. Iterating this with the different disjoint supports of the weight-$2$ codewords of $\Cc$ we get the decomposition in \eqref{eq: decomposition in projectives}. An indecomposable self-dual code of length $n \geq 4$ cannot have a weight-$2$ codeword because it would split. Therefore, for $n\geq 4$, indecomposable implies projective.

 Denote by $\Pp_{n,q}$ the number of projective self-dual codes of length $n$ over $\F_q$ and by $g_2 := \Gg_{2,q}$ the number of self-dual codes of length $2$. We set $\Pp_{0,q}=1$. We have that \begin{equation*}
    g_2 = \begin{cases}
    1 & \text{if } q \text{ is even},\\
    2 & \text{if } q \text{ is odd and } q\equiv 1 \,\,(\mathrm{mod} \, 4)\\
    0 & \text{if } q \text{ is odd and } q\equiv 3 \,\,(\mathrm{mod} \, 4)
    \end{cases}
\end{equation*}
The ways to pick $j$ disjoint pairs among $n$ coordinates is $\binom{n}{2j}(2j-1)(2j-3)\cdots3\cdot1 = \binom{n}{2j}\frac{(2j)!}{2^j j!} = \frac{n!}{(n-2j)!2^j j!}$, so using the decomposition in \eqref{eq: decomposition in projectives} we get
\begin{equation*}
    \Gg_{n,q}  = \sum_{j=0}^{n/2}\frac{n!}{(n-2j)!2^j j!}g_2^j\Pp_{n-2j,q}
\end{equation*} A direct consequence of the previous formula is that every self-dual code is projective over $\F_q$ if $q \equiv3 \pmod4$.
For $n\geq 4$, an indecomposable self-dual code is projective, so 
\begin{equation*}
    \Cc_{n,q} \leq \Pp_{n,q} \leq \Gg_{n,q}
\end{equation*}
Hence, we obtain
\begin{equation*}
    \frac{\Cc_{n,q}}{\Gg_{n,q}} \leq \frac{\Pp_{n,q}}{\Gg_{n,q}}\leq 1 \qquad \frac{\Cc_{n,q}}{\Gg_{n,q}} \leq \frac{\Cc_{n,q}}{\Pp_{n,q}}\leq 1
\end{equation*}
By Theorem \ref{th: asymptotic}, the left term in both equation tends to 1, so we get that both middle terms tend to $1$. We sum this result in the following theorem.
\begin{theorem}
     Let $n$ run over the even integers if $q$ is even or $q\equiv 1 \pmod 4$ and over the multiples of $4$ if $q\equiv 3 \pmod 4$. Then, 
       \begin{equation*}
        \lim_{n\to\infty} \frac{\Cc_{n,q}}{\Gg_{n,q}}= \lim_{n\to\infty} \frac{\Pp_{n,q}}{\Gg_{n,q}}= \lim_{n\to\infty} \frac{\Cc_{n,q}}{\Pp_{n,q}}=1.
        \end{equation*}
        In other words, asymptotically almost all self-dual codes are projective and almost all projective self-dual codes are indecomposable.
\end{theorem}

\begin{remark}
In the case of a self-dual code with a generator matrix having proportional columns, one should work with fat points to take into account multiplicities. Our approach partially works in this generality, in particular we can compute the dimension of the Schur square of a code regardless of the presence of proportional columns. However, to apply \cite[Theorem 7.3]{eisenbud2000projective} the scheme of points in $\PP^{k-1}(\K)$ must have degree $2k$ and this almost never happens when we have proportional columns.
\end{remark}

\section*{Acknowledgements}
The first version of this paper was posted online a month after~\cite{kreuzer2026codeequivalencepointset}, and at the time of writing we were not aware of that work. We are grateful to its author, Martin Kreuzer, for bringing it to our attention, as well as for his valuable suggestions and fruitful discussions.

\section*{Tool and computational resource disclosure}
The computations of Schur squares, vanishing ideals and Hilbert functions in
the examples were carried out with \textsc{SageMath}. The authors used the large language model Claude (Anthropic, model Fable 5.1) to find the code of
Example~\ref{ex:isodual}, whose properties were independently verified, and to check the text for grammatical errors and misprints. All text and mathematics were written, reviewed and edited by the authors, who take full responsibility for the content of the paper.

\bibliographystyle{plain}
\bibliography{biblio}
\end{document}